\documentclass[final]{siamltex}

\usepackage{amssymb,amsbsy}

\newcommand{\R}{\mathbb R}

\newcommand{\ink}{\rule{.5\baselineskip}{.55\baselineskip}}

\newtheorem{remark}{Remark} 
\newtheorem{example}{Example}

\title{An Explicit Cross Entropy Scheme for Mixtures\thanks{This research is supported in part by
the National
Science Foundation (NSF-DMS-1008331).}}

\author{Hui Wang\thanks{Division of Applied Mathematics, Brown University, Providence, Rhode Island 02912, USA ({\tt huiwang@cfm.brown.edu}).}
        \and Xiang Zhou\thanks{Department of Mathematics, City University of Hong Kong, Kowloon, Hong Kong, P.R. China ({\tt xiang.zhou@cityu.edu.hk}).}}

\begin{document}
\maketitle

\begin{abstract}
The key issue in importance sampling is the choice of the alternative sampling distribution, which is often chosen from the exponential tilt family of the underlying distribution. However, when the problem exhibits certain kind of nonconvexity, it is very likely that a single exponential change of measure will never attain asymptotic optimality and may lead to erroneous estimates. In this paper we introduce an explicit iterative scheme  which combines the traditional cross-entropy method and the EM algorithm to find an efficient alternative sampling distribution in the form of  mixtures. We also study the applications of this scheme to option price estimation.
\end{abstract}



\pagestyle{myheadings}
\thispagestyle{plain}
\markboth{H. WANG AND X.ZHOU}{CROSS-ENTROPY FOR MIXTURES}

\section{Introduction}

Importance sampling is a variance reduction technique in Monte Carlo simulation. It is particularly powerful in the context of rare event simulation \cite{asm3, blagly, dupwan6, hei,  sad3, sha2, sie} and has numerous applications in computational finance \cite{gla, glaheisha1, guarob, wan}. The main idea of importance sampling is to generate samples from an alternative sampling distribution and multiply the outcome with the appropriate likelihood ratio in order to obtain unbiased estimates. To be more concrete, consider a generic problem of estimating 
$$
\mu = E[V(X)],
$$
where $X$ is an $\R^d$-valued random variable defined on some probability space $(\Omega,  \mathbb{F},P)$ and $V:\R^d\rightarrow \R$ is some given nonnegative function. For simplicity, further assume that $X$ has a density function $f$. Observe that for any positive density function $g$,
\begin{equation}\label{eqn:unbias}
\mu = \int_{\R^d} V(x) f(x)\,dx = \int_{\R^d} V(x) \frac{f(x)}{g(x)}\, g(x)dx.
\end{equation}
Importance sampling  estimates $\mu$ by  the sample average of independent identically distributed (iid) copies of 
$$
\hat \mu = V(X)\frac{f(X)}{g(X)},
$$
where, abusing notation, $X$ is a representative sample from the new density function $g$. The estimate $\hat \mu$ is clearly unbiased, thanks to (\ref{eqn:unbias}). 

The efficiency of the importance sampling scheme depends on the alternative probability distribution $g$. A popular class of alternative probability measures, suggested by the large deviation analysis,  is the so-called ``exponentially tilted distribution", where 
$$
g(x) = e^{\langle \theta,x\rangle - H(\theta)},~~~~~~H(\theta) = \log \int_{\R^d}  e^{\langle \theta, x\rangle}f(x)\,dx,
$$ 
for some $\theta\in\R^d$. The cross-entropy method, invented by Reuven Rubinstein \cite{rub1, rub2, rubkro2}, is a simulation based technique for selecting the  tilting parameter $\theta$. It is easy to implement and does not require much additional 
computational cost. The implicit assumption for the resulting scheme to work well is that the exponential tilt family contains a member which serves as a    good  sampling distribution.

However, it is well known that using a single exponentially tilted distribution will lead to suboptimality and erroneous estimates \cite{glawan, deb} except in very simplistic settings. Therefore, it is in general fruitless to look for a good alternative sampling distribution within the exponential tilt family. One has to enlarge the class to have a chance of finding a good sampling distribution.

A possible solution to a general importance sampling problem is  to  use the game/subsolution framework developed in \cite{dupwan3, dupwan5}, where the resulting alternative sampling distributions are usually state-dependent. Even though it can be shown that such schemes are   optimal in an appropriate sense, the implementation requires the user to construct suitable subsolutions to certain partial differential equations, which is itself a highly nontrivial task for many practical applications with nonhomogeneous dynamics.

In this paper we consider a different approach, where the class of alternative sampling distributions consists of mixtures of exponentially tilted distributions. In many practical applications, the class of mixtures is often sufficient due to the special structures associated with the problems under consideration \cite{buc, sadbuc}. Using cross-entropy to determine a mixture sampling distribution has been considered in \cite{brechakro, blashi}. But the approaches there are very different. In \cite{brechakro} the authors mainly consider rare event probabilities and the weights are estimated from the samples by dividing the set of interest. It is not clear whether such weights are optimal. Moreover, their approach requires samples from approximate zero-variance change of measure, which is achieved from MCMC sampling and is itself difficult. \cite{blashi} is concerned with estimating the tail probabilities of the sum of iid heavy-tailed random variables. There the components of the mixture  are given and  only the weights are optimized through the classical iterative cross-entropy scheme.  The main purpose of \cite{blashi} is not developing general schemes but rather analyzing the efficiency of mixture sampling distributions in the context of heavy-tailed distributions. Our scheme is different in that we combine cross-entropy with the classical Expectation-Maximization (EM) algorithm \cite{demlairub} to produce a good mixture sampling distribution for general models.
Our scheme leads to explicit iterations to find both the optimal weight and the optimal exponential tilting parameter  for each component of the mixture.

\vspace{0.5cm}
\noindent  {\bf Comments on notation:} Here we clarify some of the convention and  notation used throughout the paper. 
\begin{enumerate}
  \item To ease exposition, in all the discussions we assume that the underlying distribution admits a density. The extension to the general distributions is straightforward. 
 
  \item The expected value with respect to the original distribution is denoted by $E$. Given an alternative density function $g$, the notation 
  $$
E_g[V(X)]
$$
denotes the expected value when $X$ has density $g$. If the alternative density belongs to some parametrized family, say $\{f_\theta\}$, then the notation is simplified to
$$
E_\theta[V(X)] = E_{f_\theta}[V(X)],
$$ 
provided there is no ambiguity.
 \item Suppose $X$ has a density $f$ and one wishes to estimate $E[V(X)]$. The alternative sampling distribution $g$ must satisfy the ``absolute continuity" condition, that is, $V(x)f(x)=0$ whenever $g(x)=0$. Throughout the paper this condition will be implicitly imposed.
  \item  $\langle \cdot, \cdot\rangle$ denotes the usual Euclidean inner product of two vectors. If the dimension is one, then the inner product is just the ordinary product.
  \item The notation
  $$
\frac{\partial}{\partial \theta}
$$
denotes the partial derivative with respect to the argument $\theta$. If $\theta$ happens to be a vector, then it denotes the gradient with respect to $\theta$.
\item The $d\times d$ identity matrix is denoted by $I_d$.
\item Inequalities between any two vector hold in the componentwise sense.
\end{enumerate} 

\section{Basics of Cross-Entropy}

In this section we describe the main idea of the cross-entropy method. A detailed account can be found in \cite{rubkro2}. To illustrate,  consider the generic problem of estimating the expected value
$$
\mu = E[V(X)],
$$
where $V$ is a nonnegative function and $X$ has density $f$. Assume that the  alternative sampling distribution is  restricted to a parameterized family, say  $\{f_\theta\}$.
A particularly popular choice is the exponential tilt family described in the Introduction.

It is well known that there exists a zero-variance alternative sampling distribution $f^*$  defined by
\begin{equation}\label{eqn:opt}
f^*(x) = \frac{1}{\mu}V(x)f(x).
\end{equation}
In other words, under $f^*$, the estimate is the constant $\mu$ and has zero variance.  Even though such a sampling distribution is impractical as it requires the knowledge of $\mu$, it leads to the heuristic principle that an alternative probability distribution  ``close" to $f^*$ should be a good choice for importance sampling.
The cross-entropy method aims to solve for the alternative probability distribution  $f_\theta$ that is {\it closest} to $f^*$ in the sense of {\it Kullback-Leibler cross-entropy} or {\it relative entropy}, defined by 
$$
R(f^*\|f_\theta) = \int f^*(x) \log \frac{f^*(x)}{f_\theta(x)}\,dx.
$$
That is, the cross-entropy method chooses  the minimizing probability density to the  minimization problem 
\begin{equation}\label{eqn:CE}
\min_{\theta} R(f^*\|f_\theta)
\end{equation}
as the alternative sampling distribution for importance sampling. 
Plugging in formula (\ref{eqn:opt}), it follows that
\begin{eqnarray*}
R(f^*\|f_\theta) &= &  \int f^*(x)\log\frac{f^*(x)}{f(x)}dx + \int f^*(x)\log\frac{f(x)}{f_\theta(x)} dx \\
& = & R(f^*\| f)  + \frac{1}{\mu}\int  V(x) \log \frac{f(x)}{f_\theta(x)}\cdot \,f(x)dx\\
& = & R(f^*\| f)  + \frac{1}{\mu}E[V(X)\log f(X)] -\frac{1}{\mu} E [V(X) \log f_\theta(X)],
\end{eqnarray*}
Since  the first term, the second term, and $\mu$ are all independent of $\theta$, the minimization problem (\ref{eqn:CE}) is equivalent to the maximization problem
\begin{equation}\label{eqn:CE_max}
 \max_\theta E [V(X) \log f_\theta(X)].
\end{equation}
The maximizing $\theta$, under mild regularity conditions, solves the equation
$$
0 = \frac{\partial}{\partial \theta} E\left[V(X)\log f_\theta(X)\right] = E\left[V(X)\frac{\partial}{\partial \theta}\log f_\theta(X) \right] .
$$
Even though there is no explicit solution to the equation in general, the cross-entropy method produces a simple  simulation based algorithm to approximate the solution. 
To be more precise,  if one replaces the expected value by sample average and  considers the corresponding stochastic version 
\begin{equation}\label{eqn:CE-basic}
0=\frac{1}{N} \sum_{k=1}^N V(X_k)\frac{\partial}{\partial \theta}\log f_\theta(X_k),
\end{equation}
where  $\{X_k\}$ are iid copies of $X$, then it is often possible to obtain a solution in {\it closed form}; see \cite[Appendix A]{rubkro2}. 
To summarize, the basic  cross-entropy method generates $N$ iid pilot samples $X_1$, $\ldots$, $X_N$ from the original distribution  and compute $\hat\theta$ from equation (\ref{eqn:CE-basic}).  Once $\hat \theta$ is obtained,  $\mu = E[V(X)]$ is estimated by importance sampling using  the alternative sampling  density  $f_{\hat\theta}$.   

 The basic cross-entropy algorithm can easily extends to a more general, iterative procedure for solving the minimization problem (\ref{eqn:CE_max}).
Every iteration in the general scheme involves two phases. In the $(i+1)$-th iteration, one first generates iid pilot  samples  from density $f_\theta$ with $\theta=\hat\theta^i$ being the current candidate of the tilting parameter.  The tilting parameter $\hat\theta^i$ is then updated to $\hat\theta^{i+1}$ based on these samples. As in the basic algorithm,  $\hat\theta^{i+1}$ often admits closed form formulas.  
If $\hat\theta^*$ is  the tilting parameter from the final iteration, then $f_{\hat\theta^*}$ is used as the alternative  sampling density to estimate $\mu$.

The iterative cross-entropy algorithm is based on the following observation. 
  Fixing an {\it arbitrary} tilting parameter say $\theta^0$,  we can rewrite the expected value in (\ref{eqn:CE_max}) as
\begin{equation}\label{eqn:CE0}
  E_{\theta^0}\left[V(X)\ell_{\theta^0}(X)\log f_\theta(X)\right],
\end{equation}
where for any $\theta$,
$$
\ell_\theta(x) = \frac{f(x)}{f_\theta(x)}
$$
denotes the likelihood ratio. 
 Consequently, the maximizer to  (\ref{eqn:CE_max}) satisfies the equation
$$
0 = E_{\theta^0}\left[V(X) \ell_{\theta^0}(X)\frac{\partial}{\partial \theta}\log f_\theta(X)\right].
$$
As before, we replace the expected value by sample average and solve the equation
$$
0= \frac{1}{N}\sum_{k=1}^N V(X_k) \ell_{\theta^0}(X_k)\frac{\partial}{\partial \theta}\log f_\theta(X_k),
$$
where $X_1,\ldots,X_N$ are iid pilot samples from  the probability density $f_{\theta^0}$. This leads to the following updating rule for $\hat\theta$.

\vspace{0.2cm}

\begin{quote}
 {\bf The updating rule of $\boldsymbol{\hat\theta}$.} Suppose that $\hat\theta^i$ is the value of the tilting parameter at the end of the $i$-th iteration.  Then $\hat\theta^{i+1}$ is set to be the solution to the equation
\begin{equation}\label{eqn:CE1}
0= \frac{1}{N}\sum_{k=1}^N V(X_k) \ell_{\hat\theta^i}(X_k)\frac{\partial}{\partial \theta}\log f_\theta(X_k),
\end{equation}
where $X_1,\ldots,X_N$ are iid pilot samples from the probability density $f_{\hat\theta^i}$.
 \end{quote}

\vspace{0.2cm} 
\noindent Equation (\ref{eqn:CE1}) is of exactly the same form as  the basic cross-entropy equation (\ref{eqn:CE-basic}),  except that  $V(X_k)$ is replaced by $V(X_k)\ell_{\hat\theta^i}(X_k)$, thus  often explicitly solvable as well.
 In particular, we have the following result concerning the normal distributions, in which case an exponentially tilted distribution  amounts to a shift in the mean. Its proof  is trivial and thus omitted.
\vspace{0.2cm}

\begin{lemma}\label{lemma:CE_normal_update} 
Suppose that  $f$ is the density of $N(0,I_d)$ and $f_\theta$ denotes the density of $N(\theta,I_d)$ for $\theta \in\R^d$. Then the solution to (\ref{eqn:CE1}) is simply
$$
\hat \theta^{i+1} = \frac{\sum_{k=1}^N V(X_k)e^{-\langle \hat\theta^i, X_k\rangle } X_k}{\sum_{k=1}^N V(X_k)e^{-\langle\hat\theta^i, X_k\rangle}} ,
$$
where $X_1,\ldots,X_N$ are iid pilot samples from   $N(\hat \theta^i,I_d)$.
  
\end{lemma}
\vspace{0.2cm}

\begin{remark}\label{remark:ini}
{\rm
Detailed discussion of the initialization of $\hat\theta^0$ can be found in \cite[Chapter 8]{rubkro2}. We will also touch upon it later in the paper. It is worth pointing out that  if the initial tilting parameter $\hat\theta^0$ is chosen such that $f_{\hat\theta^0} = f$,
then the first step of iterative cross entropy algorithm is exactly the basic cross entropy algorithm since $\ell_{\hat \theta^0}(x) =1$ for all $x$.
}
\end{remark}

\section{Cross Entropy for Mixture}\label{section:CE-mix}

It is often the case that a single exponential change of measure will not be sufficient for constructing efficient importance sampling algorithms. This was first observed in \cite{glawan} and has motivated  extensive investigations afterwards, which have led to schemes such as dynamic importance sampling \cite{dupwan3, dupwan5}  and Lyapunov function method \cite{blagly}. As powerful as these methods are, they are more accessible when the underlying dynamics are piecewise homogeneous (such as queueing models), for it is easier to explore the structure of some related partial differential equations. In financial models, however, due to the nonlinearity commonly associated with diffusive dynamics, it is very difficult to extract information from such partial differential equations.

In this section, we consider an alternative method based on mixtures, which take the following form. Let $\{f_\alpha(x): x\in\R^d\}$ be a family of density functions.  Given a positive integer $m$, let $\theta = (w_1,\ldots, w_m; \alpha_1,\ldots,\alpha_m)$, where $(w_1,\ldots,w_m)$ is a vector of weights, that is,
$$
w_i> 0,~~~ w_1+\cdots+w_m=1,
$$
and  define the mixture density 
$$
h_\theta(x) = \sum_{j=1}^m w_j \cdot f_{\alpha_j}(x).
$$ 
In practice, $\{f_\alpha\}$ is usually chosen to be the family of  exponentially tilted distributions, and $\{h_\theta\}$ represents their mixtures. 

Recall that we are interested in estimating the expected value  $\mu = E[V(X)]$, 
where $V$ is nonnegative and $X$ has density $f$. Fix an arbitrary parameter $\theta^0  = (w_1^0,\ldots, w^0_m; \alpha^0_1,\ldots, \alpha^0_m)$. As before, the  cross-entropy scheme amounts to solving the maximization problem (\ref{eqn:CE0}), namely,
$$
\max_\theta E_{\theta^0}[V(X)\ell_{\theta^0}(X)\log h_\theta(X)],
$$
where for any $\theta$
$$
\ell_\theta(x) \triangleq  \frac{f(x)}{h_\theta(x)}.
$$
The stochastic version of the previous maximization problem is to maximize 
$$
\frac{1}{N}\sum_{k=1}^N \bar V(X_k)\log h_\theta(X_k),
$$
where $X_1,\ldots,X_N$ are iid samples from the mixture  density $h_{\theta^0}(x)$ and 
\begin{equation}\label{eqn:defofbarV}
\bar V(x) \triangleq V(x)\ell_{\theta^0}(x).
\end{equation}
The problem here is that  this maximization problem will not have explicit solution anymore.

To resolve this issue,  we resort to the idea of EM algorithm. Given any $\theta$,  define (abusing  notation) the joint probability distribution function by
$$
h_\theta(x,j) \triangleq w_j \cdot f_{\alpha_j}(x), ~~~\mbox{for $x\in\R^d$ and $j=1,\ldots,m$}.
 $$
Note that the marginal of $h_\theta(x,\cdot)$ on $x$ is exactly the mixture $h_\theta(x)$. Now we can write, for any $j=1,\ldots,m$, 
\begin{eqnarray*}
\bar V(X_k)\log h_\theta(X_k) &= & \bar V(X_k) \log  h_\theta(X_k,j) - \bar V(X_k) \log\frac{h_\theta(X_k,j)}{h_\theta(X_k)}\\
& = &\bar  V(X_k) \log h_\theta(X_k,j) -\bar V(X_k) \log  h_\theta(j|X_k),
\end{eqnarray*}
where  $h_\theta(j|x) =h_\theta(x,j)/h_\theta(x)$ is the conditional probability distribution function. 
Since the above display is true for any possible value of $j$, we can integrate both sides against the conditional distribution  $h_{\theta^0}(\cdot|x)$ and obtain
\begin{eqnarray*}
\bar V(X_k)\log h_\theta(X_k)& =& \sum_{j=1}^m \bar  V(X_k) \log h_\theta(X_k,j) \cdot h_{\theta^0}(j|X_k) \\
&& ~~~~ -\sum_{j=1}^m \bar V(X_k) \log  h_\theta(j|X_k)\cdot h_{\theta^0}(j|X_k) \\
& =& \sum_{j=1}^m \bar  V(X_k) \log  h_\theta(X_k,j)  \cdot h_{\theta^0}(j|X_k) \\
&& ~~~~~ - \sum_{j=1}^m \bar V(X_k) \log  h_{\theta^0}(j|X_k)\cdot h_{\theta^0}(j|X_k)  \\
& & ~~~~~ +\sum_{j=1}^m \bar V(X_k) \log \frac{ h_{\theta^0}(j|X_k)}{h_{\theta}(j|X_k)}\cdot h_{\theta^0}(j|X_k).
\end{eqnarray*}
Following the ideas of EM algorithm, we choose $\theta$ that maximizes the summation (over $k$) of the first term  on the right-hand-side, since  the second  term is independent of $\theta$ and the last term is nonnegative (which is essentially a relative entropy itself and takes value zero if $\theta=\theta^0$). That is, we choose $\theta$ to maximize
\begin{eqnarray*}
\lefteqn{\sum_{k=1}^N \sum_{j=1}^m \bar  V(X_k) {\mathbf  {log} !!} h_\theta(X_k,j) \cdot h_{\theta^0}(j|X_k)} \\
& = & \sum_{k=1}^N \sum_{j=1}^m \bar V(X_k) \log [w_j \cdot f_{\alpha_j}(X_k)]\cdot h_{\theta^0}(j|X_k)\\
& = &  \sum_{j=1}^m\sum_{k=1}^N  \bar V(X_k) \log f_{\alpha_j}(X_k) \cdot h_{\theta^0}(j|X_k)\\
& & ~~~~~+ \sum_{j=1}^m\sum_{k=1}^N  \bar V(X_k) \log w_j\cdot  h_{\theta^0}(j|X_k).
\end{eqnarray*}
Note that this maximization problem is completely separated. In other words, for each $j=1,\ldots,m$, we choose $\alpha_j$ to maximize
\begin{equation}\label{eqn:CE1_mix}
\sum_{k=1}^N \bar V(X_k) h_{\theta^0}(j|X_k) \log f_{\alpha_j}(X_k),
\end{equation}
which often admits explicit formula, and we choose $(w_j)$ to maximize
$$
\sum_{j=1}^m \sum_{k=1}^N  \bar V(X_k) h_{\theta^0}(j|X_k) \log w_j ,~~~~\mbox{s.t. } w_j\ge 0,~\sum_{j=1}^m w_j=1,
$$
which can be solved explicitly and the maximizer is
\begin{equation}\label{eqn:CE1_mix_weight}
w_j^* = \frac{\sum_{k=1}^N \bar V (X_k)  h_{\theta^0}(j|X_k) }{\sum_{l=1}^m\sum_{k=1}^N \bar V(X_k)   h_{\theta^0}( l|X_k)}.
\end{equation}
Recalling the definition of $\bar V$ in (\ref{eqn:defofbarV}), we arrive at  an iterative cross-entropy scheme with the following updating rule.
\vspace{0.2cm}

\begin{quote}
 {\bf The updating rule of $\boldsymbol{\hat\theta}$.} Suppose that $\hat\theta^i$ is the parameter for the mixture at the end of the $i$-th iteration.  Define
$$
h_{\hat\theta^i}(j|x) = \frac{h_{\hat\theta^i}(x,j)}{h_{\hat\theta^i}(x)}.
$$ 
Then $\hat\theta^{i+1}$ is set to be  $(\hat w_1^{i+1}, \ldots,\hat w_m^{i+1}; \hat \alpha_1^{i+1},\ldots,\hat \alpha_m^{i+1}) $, where for every $j=1,\ldots,m$,
\begin{equation}\label{eqn:CE1_w}
\hat w_j^{i+1} = \frac{\sum_{k=1}^N  V(X_k)\ell_{\hat\theta^i}(X_k)  h_{\hat \theta^i}(j|X_k) }{\sum_{l=1}^m\sum_{k=1}^N  V(X_k)\ell_{\hat\theta^i} (X_k)   h_{\hat \theta^i}( l|X_k)},
\end{equation}
and $\hat\alpha_j^{i+1}$ is the solution to the equation
\begin{equation}\label{eqn:CE1_mix_update}
0= \frac{1}{N}\sum_{k=1}^N  V(X_k)\ell_{\hat\theta^i}(X_k) h_{\hat\theta^i}(j|X_k) \frac{\partial}{\partial \alpha}\log f_{\alpha}(X_k),
\end{equation}
where $X_1,\ldots,X_N$ are iid pilot samples from the mixture density $h_{\hat\theta^i}$.
 \end{quote}

\vspace{0.2cm} 

\noindent Equation (\ref{eqn:CE1_mix_update}) is of the same form as (\ref{eqn:CE1}) and often admits explicit formula, especially when $\{f_\alpha\}$ is chosen to be the family of exponentially tilted distributions. The most interesting case for us is concerned with the normal distribution, where we have a result similar to Lemma \ref{lemma:CE_normal_update}.
\vspace{0.2cm}

\begin{lemma}\label{lemma:CE_mix_normal_update} 
Suppose that  $f$ is the density of $N(0,I_d)$ and $f_\alpha$ denotes the density of $N(\alpha,I_d)$ for $\alpha \in\R^d$. Then the solution to (\ref{eqn:CE1_mix_update}) is simply
$$
\hat \alpha_j^{i+1} = \frac{\sum_{k=1}^N V(X_k)\ell_{\hat\theta^i}(X_k) h_{\hat \theta^i}(j|X_k)\cdot X_k}{\sum_{k=1}^N V(X_k)\ell_{\hat\theta^i}(X_k) h_{\hat \theta^i}(j|X_k)} ,
$$
for each $j=1,\ldots,m$.
  
\end{lemma}
\vspace{0.2cm}

\section{Comments on Initialization}

The initialization of the cross-entropy method involves the specification of $m$ and $\hat \theta^0 = (\hat w_1^0,\ldots,\hat w_m^0; \hat \alpha_1^0,\ldots, \hat \alpha_m^0)$, that is, the number of components in the mixture, as well as  the respective weight and initial tilting parameter for each component. Note that if two initial tilting parameters (say $\hat \alpha_1^0$ and $\hat \alpha_2^0$) are set to be the same, then $\hat \alpha_1^i$ and $\hat \alpha_2^i$ will remain the same in each iteration, which essentially amounts to using fewer number of components in the mixture. Therefore, without loss of generality, we require  that $\hat \alpha_1^0,\ldots,\hat \alpha_m^0$ be all different. 

The choice of $m$ usually depends on the structure of the problem. A rule of thumb or general guidance in many financial applications is that if the function $V$ only takes positive values on the union of a collection of convex sets such as
$$
\cup_{j=1}^{m'} A_j,
$$ 
then it is often beneficial to set $m=m'$. We should illustrate this from numerical examples in the next section. Given $m$, the choice of the initial weights $(\hat w_1^0,\ldots,\hat w_m^0)$ seems to be most direct -- setting $\hat w_1^0=\cdots = \hat w_m^0=1/m$ is usually a good strategy. 

 In the following discussion, we will focus on different approaches to  initialize $(\hat\alpha_1^0,\ldots,\hat\alpha_m^0)$.  

\subsection{Initialization by  perturbation}

 In iterative cross-entropy schemes without mixture,  choice of the initial parameter can be rather flexible. A common  practice is to let it be the one that corresponds to the original distribution. An immediate extension of this strategy to schemes with mixture is to combine it with random perturbation, since we require the tilting parameters to be different. For example, if a tilting parameter $\beta$ is such that
$$
f_{\beta}=f,
$$
then one can let 
$$
\alpha^0_j = \beta + \varepsilon_j, ~~j=1,\ldots,m,
$$
where $\varepsilon_j$'s are some small random perturbations. 

For illustration, let us consider a very simple problem of estimating the probability of 
\begin{equation}\label{eqn:illus_ex}
P\{X\ge a \mbox{ or } X\le b\} = P\{X\in   [a,\infty)\cup (-\infty,b]\},
\end{equation}
where $X$ is a standard normal random variable and $b<0<a$. In this case, $f$ is the density of $N(0,1)$ and the family of exponentially tilted distributions consists of 
$$
f_\alpha = \mbox{density of }N(\alpha,1),~~~\alpha\in \R.
$$
Since the target set is the union of two convex sets, we let $m=2$ and $\hat w_1^0=\hat w_2^0 =1/2$. As for the initial tilting parameter, we simply let
$$
\hat\alpha_1^0 = 0,~~~\hat\alpha_2^0=-0.1.
$$
The value of $\hat\alpha_2^0$ is not essential. Setting $\hat\alpha_2^0$ to   other similar small values will not alter the discussion below.

Table 1 is the simulation results under various parameter settings, using $N=20000$ (pilot sample size), $\mbox{IT\_NUM}=5$ (number of cross-entropy iterations), and $n=1000000$ (sample size for importance sampling). For comparison, plain Monte Carlo simulations with the same sample size are also performed and the entry ``Var Ratio" is the ratio of the empirical variance of the plain Monte Carlo estimate to that of the cross-entropy estimate. The larger this ratio, the more significant the variance reduction from the cross-entropy algorithm.  The entries $(w_1,w_2)$ and $(\alpha_1,\alpha_2)$ are the final weights and tilting parameters from the iterative cross-entropy scheme.
\vspace{0.2cm}

{\footnotesize
\begin{center}
\begin{tabular}{|c|c|c|c|c|c||c|} \hline
 $\{a,b\}$ & True Value &    CE Est & Var Ratio & $(w_1,w_2)$ &  $(\alpha_1,\alpha_2)$ & \\ \hline\hline
  $\{1,-1.5\}$ & 0.2255 &   0.2253 & 2.4 &$(0.69, 0.31)$&    $( 1.52,  -1.72)$& {\bf A} \\ \hline
   & & 0.0289 & 13.9 & $(0.79,0.21)$ & $( 2.37,   -2.83)$ & {\bf B} \\ \cline{3-7}
  $\{2,-2.5\}$ & 0.0290 &0.0227 &  23.6 &  $(0.70,0.30)$ & $( 2.37,    2.37)$& {\bf C}\\ \cline{3-7}
  & & 0.0551 & 0.00005 & $(0.64, 0.36)$ & $( 2.37, 2.37)$&{\bf D}  \\ \hline
   & & 0.0227 & 19.7 & $(0.57,0.43)$ & $(2.37,    2.37)$ &{\bf E} \\ \cline{3-7}
  \raisebox{1.5ex}[0pt]{$\{2,-3\}$} & \raisebox{1.5ex}[0pt]{0.0241}  & 0.0475 & 0.00004 & $(0.56,0.44)$ & $(2.38,2.38)$ &{\bf F} \\ \hline
\end{tabular}
\vspace{0.2cm}

{\small Table 1. Estimating $P\{X\ge a\mbox{ or }X\le b\}$}
\end{center}
}
\vspace{0.2cm}

From this numerical experiment, one can say that this naive assignment of initial parameters works in very limited situations. For instance, if $a=1$ and $b=-1.5$  ({\bf A}), the cross-entropy algorithm consistently yields good performance. However, when $a=2$ and $b=-2.5$, the results from the cross-entropy algorithm are not consistent. Sometimes, it yields great results and the two final tilting parameters $(\alpha_1,\alpha_2)$ are where we expect them to be ({\bf B}). Sometimes, the two final tilting parameters $(\alpha_1,\alpha_2)$ are almost indistinguishable. When this happens, simulation results are either seemingly ``very accurate" with a very small variance and a large variance reduction ({\bf C}), or ``very inaccurate" with a variance much larger than that of the plain Monte Carlo estimate ({\bf D}).
This phenomenon has been observed in the literature \cite{glawan,dupwan3}. Indeed, the estimate in {\bf C} looks very accurate with the 95\% confidence interval being $0.0227 \pm 0.0001$, which significantly underestimate the true value because it does not contain any (rare) samples that fall  into the region $(-\infty,b)$. When such rare samples do emerge, the variance   rises drastically, since each such sample carries a huge likelihood ratio ({\bf D}). Similar simulation results are repeated in {\bf E} and {\bf F}, except that in this case it does not even have a good run like {\bf B}.

In summary, the assignment of initial parameters by perturbation should be used with extreme caution, especially when different components of the mixture tend to collapse.

\begin{remark} {\rm Table 1 also serves as a motivation for developing mixtures for importance sampling, since the simulation runs {\bf C}, {\bf D}, {\bf E}, and {\bf F} are essentially classical cross-entropy without using mixtures. One could easily be deceived by the ``accurate" numerical results ({\bf C} and {\bf E}).
}
\end{remark}

\begin{remark}\label{remark:rare_event}
{\rm It should be mentioned that when the estimation problem is associated with rare events, the naive initialization by perturbation will not work. For example, consider the case where the quantity of interest is a rare event probability $P\{X\in A\}$, which corresponds to $V(x) = 1_{\{x\in A\}}$. If one uses a small random perturbation to set up the initial parameters, then the denominators in the updating rule (\ref{eqn:CE1_w}) and (\ref{eqn:CE1_mix}) will most likely become zero, rendering the first iteration and hence the scheme meaningless. Therefore, in this case, other methods of initialization, such as those discussed in the next two subsections, are needed. 
}
\end{remark}

\subsection{Initialization by cross-entropy}\label{section:ini_CE}

 An observation from the perturbation method is that if the initial assignment of the  parameters is  not too far away from optimality, the cross-entropy will most likely succeed. This leads to an approach that is very similar to \cite[Chapter 8]{rubkro2}. Even though their original algorithm is developed for estimating rare event probabilities, it can be modified to  deal with mixtures. 

The idea is to introduce a {\it separate} cross-entropy iterative scheme to gradually push the parameters towards reasonably good values. Note that the {\it terminal} value of the parameters from this separate cross-entropy scheme will become the {\it initial} assignment for the cross-entropy scheme  in Section \ref{section:CE-mix}. 

To describe this scheme, we observe that in many financial applications   the function $V(x)$ is of the form
$$
V(x) = H  (x) 1_{\{F(x)\in \cup_{j=1}^m A_j\}},
$$
where $H(x)$ is  nonnegative, $F(x)$ is some function, and $A_j$'s are often associated with strike prices or barriers. Moreover, we assume that $V$ can be embedded into a collection of functions (abusing notation) indexed by $\boldsymbol{\delta}=(\delta_1,\ldots,\delta_m)$:
$$
V_{\boldsymbol{\delta}}(x) = H_{\boldsymbol{\delta}}  (x) 1_{\{F(x)\in \cup_{j=1}^m A_j(\delta_j) \}}
$$
where $H_{\boldsymbol{\delta}}(x)$ is nonnegative and   $A_j(\delta_j)$ is a sequence of decreasing sets (with respect to $\delta_j$) for each $j$ --- for this reason, ${\boldsymbol{\delta}}$ is referred to as the {\it rarity parameter}. Without loss of generality, let $V(x) = V_{\boldsymbol{\delta}}(x)$ when ${\boldsymbol{\delta}}=\boldsymbol{1} \triangleq (1,1,\ldots,1)$.

In order to avoid the collapse of the mixture, consider an increasing sequence 
$$
{\boldsymbol{\delta}}^1\le {\boldsymbol{\delta}}^2\le\cdots\le \boldsymbol{1}.
$$
If ${\boldsymbol{\delta}}^1$ is small, then  with an initial parameter assignment  from random perturbation, say $\theta^0$, the cross-entropy scheme can identify a good mixture, say $h_{\theta^1}$, for estimating
$$
E[V_{{\boldsymbol{\delta}}^1}(X)].
$$
Using $\theta^1$ as the initial input, the cross-entropy scheme can then identify a good mixture, say  $h_{\theta^2}$, for estimating
$$
E[V_{{\boldsymbol{\delta}}^2}(X)],
$$
provided ${\boldsymbol{\delta}}^2$ is not too far from ${\boldsymbol{\delta}}^1$. The iterative scheme repeats these steps until a good mixture for
$$
  E[V_{\boldsymbol{1}}(X)]=E[V(X)]
$$
is found.

The choice of the sequence $\{{\boldsymbol{\delta}}^i\}$ is to ensure that the mixture does not collapse, or that samples drawn from the mixture $h_{\theta^i}$ reach the set $A_j(\delta^{i+1}_j)$ with non-trivial probability for each $j=1,\ldots,m$. That is, we require
\begin{equation}\label{eqn:rho}
P \{F(X) \in A_j(\delta^{i+1}_j) \} \ge \rho ,~~~~j=1,\ldots,m,
\end{equation}
where $X$ has the distribution $h_{\theta^i}$ and $0<\rho<1$ is some prefixed fraction. In this case, if the pilot sample size is $N$, then the number of samples to reach $A_j(\delta^{i+1}_j)$ is approximately bounded from below
\begin{equation}\label{eqn:samples}
N \rho \cdot w,
\end{equation}
where $w$ is the corresponding weight (it is a component of $\theta^i$).

 Even though there are no explicit formula to determine $\{{\boldsymbol{\delta}}^i\}$, they can be estimated by the same pilot samples that are used in the cross-entropy iterations. First of all, in all the iterations, the weights will be fixed at $w_1=\cdots=w_m=1/m$. The reason for that is (\ref{eqn:samples}) indicates the number of samples might be too few if the weights are allowed to adapt and decrease; see Remark \ref{remark:weights_fixed}. Secondly, suppose that $X_1,\ldots,X_N$ are pilot samples from $h_{\theta^i}$, then ${\delta}^{i+1}_j$ can be approximated by the largest $\delta$ such that 
the number of samples $\{F(X_1),\ldots,F(X_N)\}$ that belong to $A_j(\delta)$ is at least $N_0$ for each $j$, where $N_0=\lfloor N\rho/m\rfloor$. Lastly, once $\boldsymbol{\delta}^{i+1} =({\delta}^{i+1}_1,\ldots, {\delta}^{i+1}_m)$ is determined,   $\theta^{i+1}$ can be obtained from equation (\ref{eqn:CE1_mix_update}) with $V(X_k)$ replaced by $
V_{\boldsymbol{\delta}^{i+1}}(X_k)$. The iteration will end if $\boldsymbol{\delta}^{i+1}\ge 1$. More precisely, we have the following pseudocode.

\begin{enumerate}
{\small
  \item[{}] {\bf Pseudocode for the initialization of  cross-entropy scheme:}
  \begin{enumerate}  \item[{}]  choose a fraction $\rho\in(0,1)$ and set $N_0=\lfloor N\rho/m\rfloor$ 
    \item[{}] set $\bar \theta^0 = (1/m,\ldots,1/m; \bar \alpha_1^0,\ldots,\bar\alpha_m^0)$ from some random perturbation
    \item[{}] set iteration counter $i=0$
    \item[(*)] generate $N$ iid pilot samples $X_1,\ldots,X_N$ from density $h_{\bar \theta^i}(x)$
    \item[{}]  set for each $j=1,\ldots,m$,\begin{equation}\label{eqn:update_delta} \bar \delta^{i+1}_j = \sup \{\delta: \#\{F(X_k) \in A_j(\delta)\}\ge N_0 \} \vee \bar\delta_j^{i}\end{equation}
    \item[{}] set $\bar{\boldsymbol{\delta}}^{i+1}= ( \bar{ {\delta}}^{i+1}_1,\ldots,\bar{ {\delta}}^{i+1}_m)$
    \item[{}] for each $j=1,\ldots,m$, set $ \bar \alpha_j^{i+1}$ as the solution to the equation
    $$
0= \frac{1}{N}\sum_{k=1}^N  V_{\bar{\boldsymbol{\delta}}^{i+1}}(X_k)\ell_{\bar\theta^i}(X_k) h_{\bar\theta^i}(j|X_k) \frac{\partial}{\partial \alpha}\log f_{\alpha}(X_k)
$$
\item[{}] set $\bar\theta^{i+1} = (1/m,\ldots,1/m; \bar \alpha_1^{i+1},\ldots,\bar\alpha_m^{i+1})$
\item[{}] set iteration counter $i=i+1$
\item[{}] if $\boldsymbol{\bar\delta}^{i}\ge 1$ then set $\hat\theta^0= \bar \theta^i$ and stop, otherwise go to (*). 

\end{enumerate}
}
\end{enumerate}

\noindent Note that in the updating equation for $\boldsymbol{\delta}$ we  enforce that its components be nondecreasing as sometimes they can be driven down by  the random noise, which might slow down the scheme. The parameter $\rho$ is usually chosen as a small number, but not too small so that   $N_0$ is of medium size (a few hundreds usually suffice).

For illustration, consider the previous numerical example of estimating the probability  in (\ref{eqn:illus_ex}). One can let $\boldsymbol{\delta}=(\delta_1,\delta_2)$ and 
$$
H_\delta(x) =1, ~~~A_1(\delta_1) =  [\delta_1 a,\infty),~~~A_2(\delta_2)= (-\infty, \delta_2 b],~~~F(x)=x.
$$
In this case, it is not difficult to check that the solution to the updating equation (\ref{eqn:update_delta}) for $\boldsymbol{\delta}$ is
\begin{eqnarray*}
\bar\delta^{i+1}_1 =\frac{X_{(N-N_0+1)}}{a} \vee \bar \delta_1^i,~~~~  \bar\delta^{i+1}_2= \frac{X_{(N_0)}}{b}\vee \bar \delta_2^i,
\end{eqnarray*}
where $
X_{(1)}\le X_{(2)}\le\cdots\le X_{(N)}
$
is the order statistics of  samples $X_1,\ldots,X_N$.  Below are the numerical results. The initial parameter is set as before, i.e.,
$$
\bar\theta^0=(0.5,0.5; 0, -0.1).
$$
We also let $\rho=0.05$. The entry $(\hat\alpha_1^0,\hat\alpha_2^0)$ in Table 2 is the {\it terminal} values of the titling parameters from the initialization scheme. Note that the weights are always (0.5,0.5) in this initialization scheme, and
$$
\hat\theta_0= (0.5,0.5; \hat\alpha_1^0,\hat\alpha_2^0)
$$
will be the initial parameter setup for the iterative cross-entropy scheme. The entry ``IT" records the number of iterations performed by the initialization scheme. We can see that the cross-entropy algorithm consistently yields accurate results with significant variance reduction.

\vspace{0.2cm}
{\footnotesize
\begin{center}
\begin{tabular}{|c|c|c||c|c|c|c|} \hline
 $\{a,b\}$ & $(\hat\alpha_1^0,\hat\alpha_2^0)$ &IT &     CE Est & Var Ratio & $(w_1,w_2)$ &  $(\alpha_1,\alpha_2)$  \\ \hline\hline
  $\{1,-1.5\}$ & $(1.13, 0.85)$ & 1 &    0.2255 & 2.4 & $(0.69, 0.31)$&    $( 1.53,  -1.87)$  \\ \hline
  $\{2,-2.5\}$ & $(2.62, -3.12)$ & 4 & 0.0290   &   13.8 &  $(0.78,0.22)$ & $( 2.38,    -2.83)$ \\ \hline
 $\{2,-3\}$& $(2.69,-3.75)$& $6$   & 0.0240 & 17.1 & $(0.95,0.05)$ & $(2.37,  -3.29)$  \\ \hline
\end{tabular}
\vspace{0.2cm}

{\small Table 2. Estimating $P\{X\ge a\mbox{ or }X\le b\}$ with CE initialization}
\end{center}
}
\vspace{0.2cm}

 The advantage of this initialization method is that it is rather mechanical and involves little analysis of the system. However, it requires some extra overhead on the computational budget. Moreover, numerical experiments seem to suggest that the terminal values of the tilting parameters from this initialization scheme sometimes are far from optimality, which slows down the convergence of the subsequent cross-entropy scheme.

\begin{remark}\label{remark:weights_fixed}
{\rm 
The weights do not have to be fixed. It is sometimes more efficient to allow the weights to adapt as long as the minimal weight remains above some prespecified threshold. 
}
\end{remark}

\subsection{Initialization by approximation}

Yet another way to set up the initial mixture distribution $\hat\theta^0$ is by approximation, which is viable in many  applications. For example, in many option pricing or risk management problems the function $V(X)$ is often nonzero only on sets of the from 
$\{R\ge a\}$ or union of such sets. The random variable $R$ can be functions of the underlying asset price, and the threshold $a$ can be some given strike price, barrier price, or loss threshold. The difficulty arises when $a$ is large. In this case, a rule of thumb is to choose a tilting parameter   in an ``economical" way so that the original distribution is tilted just enough to ensure that $\{R\ge a\}$ is barely reached with nontrivial probability. This usually means  to find an $\alpha$ 
$$
E_{f_\alpha} [R] = a.
$$
If we do this for each of such sets, then we have an initial assignment of $(\hat \alpha_1^0,\ldots,\hat\alpha_m^0)$.

Again, we use problem (\ref{eqn:illus_ex}) to illustrate this approach. In this case, what the initialization amounts to is to shift the mean of the new distribution to $a$ or $b$. That is,
$$
\hat\theta^0 = (0.5, 0.5; a, b).
$$
The results are similar to that of Table 2.
\vspace{0.2cm}

{\footnotesize
\begin{center}
\begin{tabular}{|c|c|c|c|c|} \hline
 $\{a,b\}$ &       CE Est & Var Ratio & $(w_1,w_2)$ &  $(\alpha_1,\alpha_2)$  \\ \hline\hline
  $\{1,-1.5\}$ &      0.2258 & 2.3 & $(0.70, 0.30)$&    $( 1.51,  -1.92)$  \\ \hline
  $\{2,-2.5\}$ &   0.0290   &   14.0 &  $(0.78,0.22)$ & $( 2.37,    -2.82)$ \\ \hline
 $\{2,-3\}$   & 0.0241 & 17.2 & $(0.95,0.05)$ & $(2.37,  -3.27)$  \\ \hline
\end{tabular}
\vspace{0.2cm}

{\small Table 3. Estimating $P\{X\ge a\mbox{ or }X\le b\}$}
\end{center}
}
\vspace{0.2cm}

The initialization by approximation requires some analysis of the system, which may be hard to carry out at times. However, if such an approximation is available, the convergence of the corresponding iterative  cross-entropy is usually very fast.

\section{Numerical Examples}

In this section we include a number of numerical experiments. In all the simulations, the pilot sample size is $N=10000$, the sample size for importance sampling is $n=100000$, and number of iterations in the cross-entropy algorithm is set to be $\mbox{IT\_NUM}=5$, unless otherwise specified. In all the tables, ``CE Est" stands for the estimate from importance sampling using the cross-entropy method,  ``R.E." denotes the empirical relative error, which is defined to be
$$
\mbox{empirical relative error} = \frac{\mbox{standard error of  estimate}}{\mbox{estimate}},
$$
and ``Var Ratio" is again the ratio of the variance of the plain Monte Carlo estimate to that of the cross-entropy estimate.

\begin{example}{\rm We consider a very simple example of estimating the price of an average price call option. More precisely, let the stock price be 
$$
\frac{dS_t}{S_t} = r\, dt + \sigma\,dW_t 
$$
under the risk-neutral probability measure, where $r$ is the risk-free interest rate and $W$ a standard Brownian motion. Consider a discretely monitored average price call option with payoff  $(\bar S-K)^+$ and maturity $T$, where $\bar S$ is the arithmetic mean
$$
\bar S= \frac{1}{d}\sum_{i=1}^d S_{t_i}
$$
for a given set of dates $0<t_1<\cdots<t_d=T$.
}
\end{example}

Since the option payoff is positive on a simple set $\{\bar S\ge K\}$, we will use importance sampling with the classical cross-entropy method, that is, without mixtures. Define for $1\le i\le d$, $Z_i= (W_{t_i}-W_{t_{i-1}})/\sqrt{t_t-t_{i-1}}$. Then $X=(Z_1,\ldots, Z_d)$ is a $d$-dimensional standard normal random vector. The stock price at time $t_i$ can be written as a function of $X$:
$$
S_{t_i} = S_0 \exp\left\{\left(r-\frac{1}{2}\sigma^2\right) t_i + \sigma \sum_{j=1}^i \sqrt{t_j-t_{j-1}} Z_j\right\}.
$$
Therefore, the discounted option payoff $e^{-rT}(\bar S-K)^+$ is a function of $X$ as well, say $V(X)$, and the option price is $E[V(X)]$.

Since the updating rule, or the solution to equation (\ref{eqn:CE1}), is explicitly given by Lemma \ref{lemma:CE_normal_update}, the only remaining question is the initial parameter $\hat\theta^0=(\hat\theta^0_1,\ldots,\hat\theta_d^0) \in \R^d$. If the strike price $K$ is not large, one can just let $\hat\theta^0=(0,\ldots,0)$,that is, starting with the original distribution. However, if $K$ increases, this initialization will not be adequate because fewer and fewer sample path will reach the strike price $K$; see Remark \ref{remark:rare_event}. In this case, we will need a different initialization. One may use the general initialization method by cross-entropy, which is sufficient and rather straightforward. Here we  illustrate how to use the method of approximation. What we would like is that $\bar S$ will exceed $K$ with nontrivial probability if $X=(Z_1,\ldots,Z_d)$ is distributed as $N(\hat\theta^0,I_d)$. This leads to equation
$$
K = E[\bar S] = \frac{1}{d}\sum_{i=1}^dS_0 \exp\left\{\left(r-\frac{1}{2}\sigma^2\right) t_i + \sigma \sum_{j=1}^i \sqrt{t_j-t_{j-1}} \hat\theta^0_j\right\}.
$$
A convenient choice is to let $\hat\theta^0=(a,a,\ldots,a)$, where $a$ is the solution to the equation
\begin{equation}\label{eqn:ini_a}
K =  \frac{1}{d}\sum_{i=1}^dS_0 \exp\left\{\left(r-\frac{1}{2}\sigma^2\right) t_i + \sigma \sum_{j=1}^i \sqrt{t_j-t_{j-1}}\cdot a\right\}.
\end{equation}
It can be easily solved by numerical methods. 

Some numerical results are presented in Table 4.  The model parameters are given by
$$
S_0= 50,~~~ r=0.05,~~~\sigma=0.3,~~~T=1, ~~~d=30,~~~t_i=\frac{i}{d}.
$$
The initial tilting parameter are obtained by solving equation (\ref{eqn:ini_a})  with the bisection method.

\vspace{0.2cm}
{\footnotesize
\begin{center}
\begin{tabular}{|c|c|c|c|c|c|} \hline
 Strike &       $K=50$ & $K=60$ & $K=70$ &  $K=80$ & $K=90$  \\ \hline\hline
 CE Est&      4.0766 & 1.0179 & 0.1917&  0.0309 &  0.0045   \\ \hline
  R.E. &   0.16\%   &   0.23\% &  0.29\% &  0.34\% & 0.38\% \\ \hline
 Var Ratio   & 9.5 & 18.7 & 58.3 &  277.9 & 1119.1 \\ \hline
\end{tabular}
\vspace{0.2cm}

{\small Table 4. Average price call option}
\end{center}
}
\vspace{0.2cm}

\noindent Clearly the variance reduction is more significant when the option becomes more out-of-the-money. This is a common phenomenon in  importance sampling, which is most efficient in the dealing with rare events. \hfill $\ink$

\vspace{0.3cm}

\begin{example}{\rm Consider a rainbow option with $d$-underlying assets, whose prices are modeled by geometric Brownian motions under the risk-neutral probability measure:
$$
\frac{dS^{(j)}_t}{S^{(j)}_t} = r\,dt + \sigma_j \,dW_t^{(j)},~~~~j=1,\ldots,d.
$$
Here $r$ is the risk-free interest rate, and $W=(W^{(1)},\ldots,W^{(d)})$ is a $d$-dimensional Brownian motion with covariance matrix $\Sigma=[\Sigma_{ij}]$ such that $\Sigma_{jj} = 1$ for  $j=1,\ldots,d$. We are interested in estimating the prices of an outperformance option with maturity $T$ and payoff
$$
\left(\max\left\{S_T^{(1)},\ldots,S_T^{(d)}\right\}-K\right)^+.
$$ 
}
\end{example}

Let $C$ be a Cholesky factorization of $\Sigma$, that is, $C$ is a lower triangular matrix such that $CC'=\Sigma$. We can write
\begin{equation}\label{eqn:Cholesky}
\frac{1}{\sqrt{T}}W_T = \frac{1}{\sqrt{T}}(W_T^{(1)},\ldots,W_T^{(d)}) = CX
\end{equation}
for some $d$-dimensional standard normal random vector $X$. The discounted option payoff is just $V(X)$, where 
$$
V(x) = \left(\max_{j=1,\ldots,d} S_0^{(j)}\exp\left\{-\frac{1}{2}\sigma_j^2T + \sigma_j\sqrt{T} (Cx)_j\right\}-e^{-rT}K\right)^+
$$
for $x\in\R^d$. 

Since the option has strictly positive payoff when either of the $d$ stock prices at maturity exceeds the strike price $K$, it is natural to use a mixture with $m=d$ components. The updating rule is explicitly given by equations (\ref{eqn:CE1_w})--(\ref{eqn:CE1_mix_update}) and Lemma \ref{lemma:CE_mix_normal_update}. Therefore, the only question remaining is the initialization. We will illustrate both methods, i.e., initialization by cross-entropy and by approximation.

\begin{enumerate}
  \item[{(i)}] {\bf Initialization by cross-entropy:} Rewrite the discounted option payoff as
  $$
V(X) = H(X) 1_{\{F(X) \in \cup_{j=1}^d A_j\}},
$$ 
where $H(x)=V(x)$ and $F(x)=x$ for $x\in\R^d$, and 
$$
A_j = \left\{x\in\R^d: S_0^{(j)}\exp\left\{\left(r-\frac{1}{2}\sigma_j^2\right)T + \sigma_j\sqrt{T} (Cx)_j\right\} >K\right\}.
$$
We can embed it into
$$
V_{\boldsymbol{\delta}}(X) = H_{\boldsymbol{\delta}}(X)1_{\{F(X) \in \cup_{j=1}^d A_j(\delta_j)\}}
$$
where 
$$
 H_{\boldsymbol{\delta}}(x) = \left(\max_{j=1,\ldots,d} \left[S_0^{(j)}\exp\left\{-\frac{1}{2}\sigma_j^2T + \sigma_j\sqrt{T} (Cx)_j\right\}-e^{-rT}\delta_jK\right]\right)^+
$$ and $A_j(\delta_j)$ is defined exactly as  $A_j$  except that $K$ is replaced by $\delta_j K$. The pseudocode for the initialization is given in Section \ref{section:ini_CE}, where the updating equation (\ref{eqn:update_delta}) for $\boldsymbol{\delta}$ can be explicitly solved as follows: Let $X_1,\ldots,X_N$ be the pilot samples in the $i$-th iteration. For each $j=1,\ldots,d$, define a sequence (which is indeed the corresponding stock prices)
$$
Y^{(j)}_k  =  S_0^{(j)}\exp\left\{\left(r-\frac{1}{2}\sigma_j^2\right)T + \sigma_j\sqrt{T} (CX_k)_j\right\},~~~k=1,\ldots,N.
$$
Denote its order statistics by $Y^{(j)}_{(1)}\le \cdots \le Y^{(j)}_{(N)}$. It is not difficult to see that the solution to the updating equation (\ref{eqn:update_delta}) is
$$
\bar\delta^{i+1}_j = \frac{1}{K}   Y^{(j)}_{(N-N_0+1)} \vee \bar \delta_j^i.
$$

  \item[{(ii)}]{\bf  Initialization by approximation:} The initial assignment of the tilting parameters $(\hat\alpha_1^0,\ldots,\hat\alpha_d^0)$ is such that
$$
E_{\hat\alpha_j^0}\left[S_T^{(j)}\right] = K.
$$
That is, if $X$ is distributed as $N(\hat\alpha_j^0,I_d)$, then
$$
E\left[S_0^{(j)}\exp\left\{\left(r-\frac{1}{2}\sigma_j^2\right)T + \sigma_j\sqrt{T} (CX)_j\right\} \right] =K,
$$
or equivalently,
$$
S_0^{(j)}\exp\left\{ r T + \sigma_j\sqrt{T} (C\hat\alpha^0_j)_j\right\}   =K.
$$
There are many solutions to this equation. It is convenient to let
$$
\hat\alpha_j^0 = C^{-1}\eta_j,
$$
where  $\eta_j\in\R^d$ and all its components are zero except the $j$-th component which equals 
$$
\frac{1}{\sigma_j\sqrt{T}}\left(\log\frac{K}{S_0^{(j)}} -rT\right).
$$
\end{enumerate}

\noindent In the numerical experimentation, we include results from both methods of initialization. The entry ``INI\_CE" denotes the results using initialization by cross-entropy, while ``INI\_AP" denotes those via initialization by approximation.

The first numerical example is concerned with a two-stock outperformance option with maturity $T=1$ and 
$$
S_0^{(1)}=50,~S_0^{(2)}=45,~r=0.03,~\sigma_1=0.1,~\sigma_2=0.15,~\Sigma=\left[\begin{array}{cc} 1 & 0.2 \\0.2 & 1 \end{array}\right].
$$
The numerical results are reported in Table 5. The two initialization methods yield almost identical estimates and relative errors.

\vspace{0.2cm}
{\footnotesize 
\begin{center}
\begin{tabular}{|l|c|c||c|c||c|c|} \hline
Strike & \multicolumn{2}{c||}{$K=50$} &  \multicolumn{2}{c||}{$K=60$}  &  \multicolumn{2}{c|}{$K=70$} \\ \cline{2-7}
& INI\_CE & INI\_AP & INI\_CE & INI\_AP & INI\_CE & INI\_AP \\ \hline
CE Est &  3.5898 & 3.5825 & 0.2768 & 0.2763 & 0.0093 &   0.0093  \\ \hline
R.E. &0.14\%  &  0.16\% & 0.27\%  &  0.28\%& 0.37\%  & 0.37\%\\ \hline
Var Ratio &  6.2 &5.2 & 27.2 & 26.9 & 445.8 & 456.5\\ \hline
\end{tabular}

\vspace{0.2cm}

{\small Table 5: Outperformance option with two underlying assets} 
 \end{center}
}
\vspace{0.2cm}

The second numerical example has $d=4$ underlying assets with the following parameters.
$$
S_0^{(1)} = 45, ~~S_0^{(2)} = 50,~~S_0^{(3)}=47,~~S_0^{(4)}=50,~~r=0.02,~~T=0.5,
$$
$$
\sigma_1=\sigma_2=0.1,~~\sigma_3=\sigma_4=0.2,~~\Sigma = \left[\begin{array}{rrrr} 1.0 & 0.3 & -0.2 & 0.4 \\ 0.3 & 1.0 & -0.3 & 0.1 \\ -0.2 & -0.3 & 1.0 & 0.5 \\ 0.4 & 0.1 & 0.5 & 1.0 \end{array}\right].
$$
We also let $\mbox{IT\_NUM}=10$ in the simulation because of the slow convergence when initialization by cross-entropy is adopted. The numerical results are reported in Table 6.

\vspace{0.2cm}
{\footnotesize 
\begin{center}
\begin{tabular}{|l|c|c||c|c||c|c|} \hline
Strike & \multicolumn{2}{c||}{$K=50$} &  \multicolumn{2}{c||}{$K=60$}  &  \multicolumn{2}{c|}{$K=70$} \\ \cline{2-7}
& INI\_CE & INI\_AP & INI\_CE & INI\_AP & INI\_CE & INI\_AP \\ \hline
CE Est &  4.6841 & 4.6722 & 0.5271 & 0.5284 & 0.0360 &   0.0362  \\ \hline
R.E. &0.17\%  &  0.17\% & 0.30\%  &  0.26\%& 0.34\%  & 0.34\%\\ \hline
Var Ratio &  3.5 & 3.5 & 15.1 & 19.9 & 143.0 & 143.6\\ \hline
\end{tabular}

\vspace{0.2cm}

{\small Table 6: Outperformance option with four underlying assets} 
 \end{center}
}
\vspace{0.2cm}

\noindent Even though both methods lead to significant variance reduction, more extensive numerical experiments appear to indicate that the convergence of the iterative cross-entropy algorithm is much faster when combined the method of initialization by approximation. This is because initialization by approximation usually gives a set of parameters closer to the optimal one.

\begin{example}
{\rm The setup is the same as the previous example. We are interested in estimating the price of a pyramid rainbow option with maturity $T$ and payoff
$$
\left( |S_T^{(1)}-K_1 | + \cdots +  |S_T^{(d)}-K_d|-K\right)^+
$$ 
where $K_1,\ldots,K_d$ and $K$ are all positive constants.
}
\end{example}

It is not difficult to see that the payoff of the option is positive on each of the following $2^d$ sets:
$$
  \left\{\pm[S_T^{(1)}-K_1]\ge 0,\cdots,\pm [S_T^{(d)}-K_d]\ge 0: \sum_{j=1}^d\pm[ S_T^{(j)}- K_j]\ge  K\right\}.
$$
In other words, in each set and for each $j=1,\ldots,d$, it is either $S_T^{(j)}\ge  K_j$ or $S_T^{(j)}\le K_j$, and the  absolute values are expressed accordingly. Therefore, it is natural to select $m=2^d$. As for initialization, one can use the method   by cross-entropy as before, but the convergence is slow. Instead, we will use the method by approximation, which is fairly easy to carry out and the convergence is much faster. 
To this end, recall the Cholesky factorization (\ref{eqn:Cholesky}) and  rewrite
$$
S_T^{(j)} = S_0^{(j)} \exp \left\{\left(r-\frac{1}{2}\sigma_j^2\right)T+\sigma_j\sqrt{T}(CX)_j\right\},
$$
where $X$ is a standard $d$-dimensional normal random vector. Consider a typical set, say
$$
\left\{S_T^{(1)} \ge K_1,\ldots, S_T^{(d)}\ge K_d: \sum_{j=1}^d [ S_T^{(j)}- K_j]\ge K\right\}.
$$ 
For this set, we would like to construct a tilting parameter $\hat\alpha\in \R^d$ such that if $X$ is distributed as $N(\hat\alpha,I_d)$, then the above set is not a rare event. Note that under the original distribution,
$$
E[S_T^{(j)}] = S_0^{(j)} \exp\{rT\},
$$
while under the new distribution 
\begin{equation}\label{eqn:new_mean}
E_{\hat\alpha}[S_T^{(j)}] = S_0^{(j)} \exp\left\{rT + \sigma_j\sqrt{T}(C\hat\alpha)_j\right\}:=x_j.
\end{equation}
The idea is to push $x_j$ beyond $K_j$ if $E[S_T^{(j)}]<K_j$, which means we would require 
$$
x_j \ge \max \left\{K_j,  S_0^{(j)} \exp\{rT\}\right\}:= \bar K_j.
$$
Obviously we should also require
\begin{equation}\label{eqn:mean}
\sum_{j=1}^d (x_j-K_j) \ge K.
\end{equation}
A simple choice is to let $x_j = \bar K_j +a$ for some  constant $a\ge 0$. Then (\ref{eqn:mean}) becomes
$$
a \ge \frac{1}{d}\left[ K-\sum_{j=1}^d (\bar K_j-K_j)\right].
$$
Therefore, choosing the minimal nonnegative $a$ that satisfies this requirement, we arrive at 
$$
x_j = \bar K_j + \frac{1}{d} \left[K-\sum_{i=1}^d (\bar K_i-K_i)\right]^+.
$$
The   formula of $x_j$ for a general set   can be obtained in exactly the same way (we omit the details) and is indeed
$$
x_j =  \bar K_j \pm \frac{1}{d} \left[K-\sum_{i=1}^d |\bar K_i-K_i|\right]^+,
$$
taking the ``$+$" or the ``$-$" sign depending on whether the constraint in the corresponding set is $S_T^{(j)}\ge K_j$ or $S_T^{(j)}\le K_j$.
Now plugging $x_j$ into equation (\ref{eqn:new_mean}), we have $\hat \alpha = C^{-1}\eta$,
where $\eta=(\eta_1,\ldots,\eta_d)'$ and
$$
\eta_j = \frac{1}{\sigma_j\sqrt{T}}\left(\log \frac{x_j}{S_0^{(j)}}-rT\right),~~~j=1,\ldots,d.
$$

Below are some numerical results. The first example has $d=2$ underlying assets. Thus we use mixtures with $m=2^d=4$ components. The parameters are
$$
S_0^{(1)}=50,~~S_0^{(2)}=45,~~K_1=55,~~K_2=50,~~r=0.03, ~~T=1,
$$
$$
\sigma_1=0.2,~~\sigma_2=0.25, ~~\Sigma=\left[\begin{array}{cc}1 & 0.3 \\ 0.3 & 1\end{array}\right].
$$
The simulation results are reported in Table 7.

\vspace{0.2cm}
{\footnotesize
\begin{center}
\begin{tabular}{|c|c|c|c|c|c|} \hline
 Strike &       $K=10$ & $K=20$ & $K=30$ &  $K=40$ & $K=50$  \\ \hline\hline
 CE Est&      9.3417 & 3.4025 & 0.9050&  0.1930 &  0.047   \\ \hline
  R.E. &   0.16\%   &   0.23\% &  0.29\% &  0.33\% & 0.36\% \\ \hline
 Var Ratio   & 3.4 & 6.3 & 16.5 &  64.6 & 262.3 \\ \hline
\end{tabular}
\vspace{0.2cm}

{\small Table 7. Pyramid option with two underlying assets}
\end{center}
}
\vspace{0.2cm}

The second example considers a model with $d=4$ underlying assets. The number of components in the mixture is $m=2^d=16$. The parameters are
$$
S_0^{(1)}=50, ~~S_0^{(2)}=45,~~S_0^{(3)}=45,~~S_0^{(4)} = 30,~~T=1,
$$
$$
K_1=55,~~K_2 = K_3=50,~~K_4=35,~~r=0.03,~~
 \sigma_1=\sigma_2=0.15,
$$
$$
\sigma_3=\sigma_4=0.2,~~\Sigma=\left[\begin{array}{cccc}1 & 0.1 & -0.2 & 0.3 \\ 0.1 & 1 & -0.5 & 0.4 \\ -0.2 & -0.5 & 1 & 0.2 \\ 0.3 & 0.4 & 0.2 & 1\end{array}\right].
$$
The simulation results are reported in Table 8.
\vspace{0.2cm}

{\footnotesize
\begin{center}
\begin{tabular}{|c|c|c|c|c|c|} \hline
 Strike &       $K=20$ & $K=30$ & $K=40$ &  $K=50$ & $K=60$  \\ \hline\hline
 CE Est&      8.8209 & 3.2507 & 0.8504 &  0.1713 &  0.032   \\ \hline
  R.E. &   0.16\%   &   0.23\% &  0.30\% &  0.35\% & 0.38\% \\ \hline
 Var Ratio   & 4.0 & 6.2 & 14.8 &  51.6 & 232.8 \\ \hline
\end{tabular}
\vspace{0.2cm}

{\small Table 8. Pyramid option with four underlying assets}
\end{center}
}
\vspace{0.2cm}

\begin{example}
{\rm 
Consider a two dimensional CEV model where the stock prices satisfy the following stochastic differential equations under the risk-neutral probability measure:
\begin{eqnarray*}
dS_t & = & r S_t  \,dt + \sigma_1 S_t^{\gamma_1}\,dW_t,\\
dH_t & = & rH_t\,dt + \sigma_2H_t^{\gamma_2}\,dB_t,
\end{eqnarray*}
where $\gamma_1,\gamma_2\in[0.5,1)$ are given constants, and $(W,B)$ is a two-dimensional Brownian motion with covariance matrix
$$
\Sigma = \left[\begin{array}{cc}1 & \rho \\ \rho & 1\end{array}\right].
$$
We are interested in estimating the price of a digital call option with payoff
$$
1_{\{\max\{c_1S_t, ~c_2H_t\}\ge K\}}
$$
with maturity $T$. 
}
\end{example}

We work with  the discounted stock prices, which makes the  future analysis slightly easier. More precisely, let $X_t = e^{-rt}S_t$ and $Y_t = e^{-rt} H_t$. It follows from It\^o formula that
\begin{eqnarray} \label{eqn:SDE_X}
dX_t & = &  \sigma_1 e^{-r(1-\gamma_1)t}X_t^{\gamma_1}\,dW_t,\\
dY_t & = &  \sigma_2 e^{-r(1-\gamma_2)t}Y_t^{\gamma_2}\,dB_t. \label{eqn:SDE_Y}
\end{eqnarray}
Since there are two ways to have a positive payoff, i.e., either of the two stock prices reaching a high enough level at maturity, we let $m=2$. The question is how to determine the initial tilting parameters. One can use initialization  by cross-entropy, but the convergence is very slow. So here we  use the method of approximation again. 

 Consider a simple change of measure $Q_1$ under which $W$ becomes a Brownian motion with constant drift $x$. 
The drift $x$ will be determined so that the equation
$$
E_{Q_1}[X_T] = e^{-rT} E_{Q_1}[S_T] = e^{-rT}K/c_1
$$
is satisfied approximately, which ensures $c_1S_T$ to  reach the target level $K$with nontrivial probability.  To this end, we rewrite equation (\ref{eqn:SDE_X}) as
$$
dX_t = x\sigma_1 e^{-r(1-\gamma_1)t}X_t^{\gamma_1}\,dt+ \sigma_1 e^{-r(1-\gamma_1)t}X_t^{\gamma_1}\,d(W_t-xt),
$$
where $ W_t -xt$ is a Brownian motion under probability measure $Q_1$. Taking expected values on both sides and using H\"older inequality, we arrive at
\begin{eqnarray*}
E_{Q_1}[X_t] &=& X_0 + \int_0^t x\sigma_1 e^{-r(1-\gamma_1)s}E[X_s^{\gamma_1}]\,ds \\
& \ge & X_0 + \int_0^t x\sigma_1 e^{-r(1-\gamma_1)s}\left(E[X_s]\right)^{\gamma_1}\,ds .
\end{eqnarray*}
This motivates the approximation of $E_{Q_1}[X_t]$ by the function $f(t)$, which is defined to be the solution to the equation
\begin{equation}\label{eqn:f}
f(t) = X_0 + \int_0^t x\sigma_1 e^{-r(1-\gamma_1)s}f^{\gamma_1}(s)\,ds.
\end{equation}
It can be shown that $f(t)$ serves as a lower bound for $E_{Q_1}[X_t]$. Actually, it is also a very good approximation for the following reason. Note that  $E[X_s^{\gamma_1}] = (E[X_s])^{\gamma_1}$ when $\gamma_1=1$ or when $X_s$ is a constant (i.e., variance of $X_s$ is zero).  Therefore, when $\gamma_1$ is close to one, we expect $f(s)$ to approximate $E_{Q_1}[X_s]$ well. When $\gamma_1$ is away from one, the variance of $X_s$ is often relatively small, which again implies that $f(s)$ is a good approximation. 
It is natural now to choose $x$ so that $f(T) = e^{-rT}K/c_1$. To this end, observe that equation (\ref{eqn:f}) is equivalent to the ordinary differential equation
$$
\frac{1}{f^{\gamma_1}} \frac{df}{dt} = x\sigma_1 e^{-r(1-\gamma_1)t},~~~f(0)=X_0=S_0.
$$
Solving this equation explicitly we obtain  
$$
f^{1-\gamma_1}(t) - S_0^{1-\gamma_1} =\frac{x\sigma_1}{r}\left[1-e^{-r(1-\gamma_1)t}\right]. 
$$
Letting $f(T) = e^{-rT}K/c_1$, we arrive at
\begin{equation}\label{eqn:x}
x = \frac{r}{\sigma_1}\frac{(e^{-rT}K/c_1)^{1-\gamma_1} - S_0^{1-\gamma_1}}{1-e^{-r(1-\gamma_1)T}}.
\end{equation}
Note that  $dB_t = \rho dW_t + \sqrt{1-\rho^2} d\bar B_t$ where $W$ and $\bar B$ are two independent standard   Brownian motions. The probability measure $Q_1$ shifts the drift of $W$ to $x$, but    $\bar B_t$ will remain a standard Brownian motion under $Q_1$.
 
We use the following version of Euler scheme to discretize the processes under $Q_1$.  Given a large integer $N_0$, let $\Delta t=T/N_0$ and $t_i=i\Delta t$ for $i=0,\ldots,N_0$. Let $\hat X_0=X_0=S_0,\hat Y_0=Y_0=H_0$, and recursively define 
\begin{eqnarray*}
\hat X_{t_{i+1}} & = &   \hat X_{t_i}+x\sigma_1 e^{-r(1-\gamma_1)t_i}\hat X_{t_i}^{\gamma_1}\Delta t+ \sigma_1 e^{-r(1-\gamma_1)t_i}\hat X_{t_i}^{\gamma_1}\sqrt{\Delta t}\cdot Z_{i+1} ,\\
\hat Y_{t_{i+1}} & = &  \hat Y_{t_i}+\rho x\sigma_2 e^{-r(1-\gamma_2)t_i}\hat Y_{t_i}^{\gamma_2}\Delta t  + \rho \sigma_2 e^{-r(1-\gamma_2)t_i}\hat Y_{t_i}^{\gamma_2}\sqrt{\Delta t}\cdot   Z_{i+1}  \\ 
& & ~~~~~~~~~+~ \sigma_2 e^{-r(1-\gamma_2)t_i}\hat Y_{t_i}^{\gamma_2}\sqrt{\Delta t}\cdot \sqrt{1-\rho^2}R_{i+1}  ,
\end{eqnarray*}
where $\{Z_i,R_i\}$ are iid standard  normal random variables.  In the implementation, if $\hat X_{t_i}$ or $\hat Y_{t_i}$ becomes negative, we reset it to zero.

The second change of measure $Q_2$ can be determined similarly. Under $Q_2$, $B$ becomes a Brownian motion with drift $y$, which,  analogous to (\ref{eqn:x}),   is given by
$$
y = \frac{r}{\sigma_2}\frac{(e^{-rT}K/c_2)^{1-\gamma_2} - S_0^{1-\gamma_2}}{1-e^{-r(1-\gamma_2)T}}, 
$$
and the corresponding discretization scheme is
\begin{eqnarray*}
\hat X_{t_{i+1}} & = &   \hat X_{t_i}+\rho y \sigma_1 e^{-r(1-\gamma_1)t_i}\hat X_{t_i}^{\gamma_1}\Delta t+ \rho \sigma_1 e^{-r(1-\gamma_1)t_i}\hat X_{t_i}^{\gamma_1}\sqrt{\Delta t}\cdot Z_{i+1} \\
& & ~~~~~~~~~+~ \sigma_1 e^{-r(1-\gamma_1)t_i}\hat X_{t_i}^{\gamma_1}\sqrt{\Delta t}\cdot \sqrt{1-\rho^2}R_{i+1}  ,\\
\hat Y_{t_{i+1}} & = &  \hat Y_{t_i}+ y\sigma_2 e^{-r(1-\gamma_2)t_i}\hat Y_{t_i}^{\gamma_2}\Delta t  +  \sigma_2 e^{-r(1-\gamma_2)t_i}\hat Y_{t_i}^{\gamma_2}\sqrt{\Delta t}\cdot   Z_{i+1} , 
\end{eqnarray*}
where $\{Z_i,R_i\}$ are iid standard normal random variables. Again, if $\hat X_{t_i}$ or $\hat Y_{t_i}$ becomes negative, we reset it to zero.

Some numerical results are reported in Table 9. The parameters are given by
$$
S_0=50,~ H_0=48,~\sigma_1=0.3,~\sigma_2 =0.35, ~\gamma_1=0.5, ~\gamma_2=0.7,~\rho=0.3,
$$
$$
c_1=c_2=1,~r=0.03, ~T=1,~N_0=50.
$$

{\footnotesize
\begin{center}
\begin{tabular}{|c|c|c|c|c|c|} \hline
 Strike &       $K=50$ & $K=55$ & $K=60$ &  $K=65$ & $K=70$  \\ \hline\hline
 CE Est&    0.8297   & 0.1908 & 0.0314  &  0.0039  & $3.3638\times10^{-4}$    \\ \hline
  R.E. &   0.14\%   &  0.37\% &  0.48\% & 0.56\% & 0.65\% \\ \hline
 Var Ratio   & 1.0  & 3.0 & 13.8 & 82.4  & 545.5 \\ \hline
\end{tabular}
\vspace{0.2cm}

{\small Table 9. Digital options for a CEV model}
\end{center}
}

%\begin{example}
%{\rm Consider the Heston stochastic volatility model where under the %risk-neutral probability measure the stock price satisfies%
%\begin{eqnarray*}
%dS_t & = & rS_t+\sqrt{\nu_t}S_t \,dW_t\\
%d\nu_t & = & a(b-\nu_t)\,dt +\sigma\sqrt{\nu_t}\,dB_t.  
%\end{eqnarray*}
%Here $a, b, \sigma$ are positive constants and $(W,B)$ is a two-dimensional %Brownian motion with covariance matrix 
%$$
%\Sigma = \left[\begin{array}{cc} 1 & \rho \\ \rho & 1\end{array}\right].
%$$
%It is assumed  that $2ab>\sigma^2$ so that the volatility $\nu$ is strictly %positive with probability one. We would like  to estimate the price of a call option with maturity $T$ and strike price $K$.
%}
%\end{example}

%\begin{enumerate}
%  \item  Other Rainbow options.
%  \item Rainbow options under the jump diffusion
%  \item Rainbow options under the CEV processes
%\end{enumerate}  

\end{document}